\documentstyle[11pt]{article}
       \font\tenmsb=msbm10
       \font\sevenmsb=msbm7
       \font\fivemsb=msbm5
       \catcode`\@=11
       \ifx\amstexloaded@\relax\catcode`\@=\active
       \endinput\else\let\amstexloaded@\relax\fi
       \def\spaces@{\space\space\space\space\space}
       \def\spaces@@{\spaces@\spaces@\spaces@\spaces@\spaces@}
       \def\space@.{\futurelet\space@\relax}
       \space@. %
       \def\Err@#1{\errhelp\defaulthelp@\errmessage{AmS-TeX error: #1}}
       \def\relaxnext@{\let\next\relax}
       \def\accentfam@{7}
       \def\noaccents@{\def\accentfam@{0}}
       \def\Cal{\relaxnext@\ifmmode\let\next\Cal@\else
       \def\next{\Err@{Use \string\Cal\space only in math mode}}\fi\next}
       \def\Cal@#1{{\Cal@@{#1}}}
       \def\Cal@@#1{\noaccents@\fam\tw@#1}
       \def\Bbb{\relaxnext@\ifmmode\let\next\Bbb@\else
       \def\next{\Err@{Use \string\Bbb\space only in math mode}}\fi\next}
       \def\Bbb@#1{{\Bbb@@{#1}}}
       \def\Bbb@@#1{\noaccents@\fam\msbfam#1}
       \newfam\msbfam
       \textfont\msbfam=\tenmsb
       \scriptfont\msbfam=\sevenmsb
       \scriptscriptfont\msbfam=\fivemsb

\def\Z{{\Bbb Z}}
\def\R{{\Bbb R}}

\def\C{{\Bbb C}}

\newtheorem{Theorem}{Theorem}
\newtheorem{Lemma}{Lemma}[section]
\newtheorem{Proposition}{Proposition}[section]

\newtheorem{Definition}{Definition}[section]

\@addtoreset{equation}{section}

\newcommand{\bq}{\begin{equation} }
\newcommand{\eq}{\end{equation} }
\newcommand{\bp}{\begin{Proposition}}
\newcommand{\ep}{\end{Proposition}}
\newcommand{\bdf}{\begin{Definition}}
\newcommand{\edf}{\end{Definition}}
\newcommand{\bl}{\begin{Lemma}}
\newcommand{\el}{\end{Lemma}}
\newcommand{\ba}{\begin{array}}
\newcommand{\ea}{\end{array}}
\newcommand{\bea}{\begin{eqnarray}}
\newcommand{\eea}{\end{eqnarray}}

\begin{document}
\setlength{\columnsep}{5pt}
\title{\bf Symplectic Energy and Lagrangian Intersection Under Legendrian
Deformations \footnotetext{MSC2000:  53D10 53D12 53D40 57R22
;}\footnotetext{ Keywords: Arnold conjecture, Floer homology,
Lagrangian intersection, Symplectic energy.}}
\author{{\bf Hai-Long Her}}
\date{}

\maketitle

\begin{quote}
\small {\bf Let $M$ be a compact symplectic manifold, and $L\subset
M$ be a closed Lagrangian submanifold which can be lifted to a
Legendrian submanifold in the contactization of $M$. For any
Legendrian deformation of $L$ satisfying some given conditions, we
get a new Lagrangian submanifold $L'$. We prove that the number of
intersection $L\cap L'$ can be estimated from below by the sum of
$\Z_2$-Betti numbers of $L$, provided they intersect transversally.}
\end{quote}

\section {Introduction.}

 \ \ \ \ \ \  In 1965, V. I. Arnold\cite{A1}\cite{A2} formulated his famous
conjectures concerning about the number of fixed points of
Hamiltonian diffeomorphisms of any compact symplectic manifold and
the number of intersection points of any Lagrangian submanifold
with its Hamiltonian deformations in a symplectic manifold. More
precisely, his conjectures can be written in topological terms as
$$\#{\rm Fix}(\psi_M)\geq\left\{
\begin{array}{ccc}
{\rm sum\ of\ Betti\ numbers\ of\ } M,&
{\rm all\ fixed\ points\ are\ nondegenerate;}\\

{\rm cuplength\ of\ } M,& {\rm some\ fixed\ points\ maybe\
degenerate,}
\end{array}
\right.
$$
and
$$\#(L\cap \psi_M(L))\geq\left\{
\begin{array}{cc}
{\rm sum\ of\ Betti\ numbers\ of\ } L,& {\rm intersection\ points\ are\ transverse;}\\

{\rm cuplength\ of\ } L,& {\rm maybe\ \ non-transverse,}
\end{array}
\right.
$$
where $M$ is a symplectic manifold, $L\subset M$ is a Lagrangian
submanifold, $\psi_M$ is a Hamiltonian diffeomorphism.

To prove these two conjectures, many works have been done, the
pioneers of them are due to Conley-Zehnder\cite{CZ}, Gromov\cite{G}
and Floer\cite{F1}-\cite{F4}. Especially, Floer originally developed
the seminal method, motivated by the variational method used by
Conley and Zehnder and the elliptic PDE techniques introduced by
Gromov, which is now called Floer homology theory, and solved many
special cases of Arnold's conjectures. In 1996, Fukaya-Ono\cite{FO}
, Liu-Tian\cite{LT} and Ruan\cite{R} independently proved the first
conjecture for general compact symplectic manifolds in the
non-degenerate case. While the conjecture for general symplectic
manifolds in the degenerate case is still open.

For the second conjecture, Floer\cite{F1}\cite{F4} gave the proof
under an additional assumption $\pi_2(M,L)=0$. We write his result
for the case that all intersections are transverse.\\ \\
\noindent {\bf Floer's Theorem}.\ \ {\it Let L be a closed
Lagrangian submanifold of a compact(or tame) symplectic manifold
$(M,\omega)$ satisfying $\pi_2(M,L)=0$, and $\psi_M$ be a
Hamiltonian diffeomorphism, then $\#(L\cap\psi_M(L))\geq {\rm
dim}H_*(L,\Z_2)$, if all intersections are transverse. }\\

In general, the condition $\pi_2(M,L)=0$ can not be removed. For
instance, let $L$ be a circle in $\R^2$, then $\pi_2(\R^2,L)\neq
0$, however, there always exists a Hamiltonian diffeomorphism
which can translate $L$ arbitrarily far from its original
position.

To prove his theorem, Floer introduced the so-called Floer
homology group for Lagrangian pairs and showed that it is
isomorphic to the homology of $L$ under the condition above. The
definition of Floer homology for Lagrangian pairs was generalized
by Oh\cite{Oh2} in the class of monotone Lagragian submanifolds
with minimal Maslov number being at least 3. However, for general
Lagrangian pairs, the Floer homology is hard to define due to the
bubbling off phenomenon and some essentially topological
obstructions \cite{FO3}, which is much different from the Hamiltonian
fixed point case.

Therefore, if we want to throw away the additional assumption, we
have to restrict the class of Hamiltonian diffeomorphisms. For the
simplest case that $\psi_M$ is $C^0$-small perturbation of the
identity, the Lagrangian intersection problem is equivalent to the
one for zero sections of cotangent bundles, which is proved by
Hofer\cite{H0} and Laudenbach-Sikorav\cite{LS}. Yu.V.
Chekanov\cite{C1}\cite{C2} also gave a version of Lagrangian
intersection theorem which used the notion of symplectic energy
introduced by Hofer\cite{H1} (for $(\R^{2n},\omega_0)$) and
Lalonde-McDuff\cite{LM} (for general symplectic manifolds).
Following their notations, we denote by ${\Cal H}(M)$ the space of
compactly supported smooth functions on $[0,1]\times M$. Any $H\in{\Cal
H}(M)$ defines a time dependent Hamiltonian flow $\phi^t_H$ on
$M$, all such time-1 maps $\{\phi^1_H,\ H\in{\Cal H}(M)\}$ form a
group, denoted by $Ham(M)$. Now we define a norm on ${\Cal H}(M)$:
$$\|H\|=\int_0^1(\max_{x\in M}H(t,x)-\min_{x\in M}H(t,x))dt,$$
and we can define the energy of a $\psi\in Ham(M)$ by
$$E(\psi)=\inf_{H}\{\|H\| \mid \psi=\phi^1_H,\ H\in {\Cal H}(M) \}.$$

For a compact symplectic manifold $(M,\omega)$, there always
exists an almost complex structure $J$ compatible with $\omega$,
so $(M,\omega,J)$ is a compact almost complex manifold, we denote by
${\Cal J}$ the set of all such $J$. Let $\sigma_S(M,J)$ and
$\sigma_D(M,L,J)$ denote the minimal area of a $J$-holomorphic
sphere in $M$ and of a $J$-holomorphic disc in $M$ with boundary
in $L$, respectively. If there is no such $J$-holomorphic curve,
these numbers will be infinity. Otherwise, minimums are obtained
by the Gromov compactness theorem\cite{G}, and they are always
positive. We write
$\sigma(M,L,J)=\min(\sigma_S(M,J),\sigma_D(M,L,J))$, and
$\sigma(M,L)=\sup_{J\in{\Cal J}}\sigma(M,L,J)$. Then Chekanov
showed the following theorem.\\ \\
\noindent{\bf Chekanov's Theorem}\cite{C2}.\ \ {\it If
$E(\psi)<\sigma(M,L)$, then $\#(L\cap \psi(L))\geq {\rm
dim}H_*(L,\Z_2)$, provided all intersections are transverse.}\\

\noindent {\it Remark}.  For the non-transverse case, under
similar assumptions, C.-G. Liu\cite{L} also got an estimate for
Lagrangian intersections by cup-length of $L$.\\

In this paper, we give an analogous Lagrangian intersection
theorem, but the Hamiltonian deformation $\psi$ will be replaced
by a ``Legendrian deformation" $\tilde{\psi}$ (which will be
explained in the sequel). In fact, K. Ono has shown such a result
still under the assumption $\pi_2(M,L)=0$.

Suppose that the symplectic structure $\omega$ is in an integral
cohomology class, and there exists a principal circle bundle
$\pi:\ N\rightarrow M$ with a connection so that the curvature
form is $\omega$, that means for a connection form $\alpha$, one
has $d\alpha=\pi^*\omega$. We see that the horizontal distribution
$\xi=Ker(\alpha)$ is a co-oriented contact structure on $N$. We
say $L$ satisfies the Bohr-Sommerfeld condition if $\alpha|_L$ is
flat, or in other words, it can be lifted to a Legendrian
submanifold
$\Lambda$ in $N$. The following is Ono's result.\\ \\
\noindent {\bf Ono's Theorem}\cite{On}.\ \ {\it Given a contact
isotopy $\{\tilde{\psi}_t\ |\ 0\le t\le 1\}$ on $N$, if $L$ is a
Lagrangian submanifold of $M$ which can be lifted to a Legendrian
submanifold $\Lambda$ in $N$, and
 $\pi_2(M,L)=0$, then $\#(L\cap\pi\circ\tilde{\psi}_1(\Lambda))\geq {\rm
dim}H_*(L,\Z_2)$, provided $L$ and
$\pi\circ\tilde{\psi}_1(\Lambda)$
intersect transversally.}\\

\noindent{\it Remark.} Since a Hamiltonian isotopy of $M$ can be
lifted to a contact isotopy of $N$, Ono's theorem is a
generalization of the previous Floer's theorem.

Eliashberg, Hofer, and Salamon\cite{EHS} also independently obtained
a result similar to Ono's theorem, they successfully overcome some
difficulties due to the non-compactness of the symplectization
manifold, while their arguments involve some complicated conditions
for avoiding bubbling off.

In the present paper, we will throw away the assumption $\pi_2(M,
L)=0$ in Ono's theorem, at the same time, we will add a certain
restrictive condition on the class of Legendrian deformation
$\tilde{\psi}$. Firstly, we denote by $\tilde{L}$ the image of
$\Lambda$ under the principal $S^1$-action on $N$. We denote by
$(SN,\omega_{\xi})$ the symplectization of the contact manifold
$(N,\xi)$ with co-oriented contact strcture $\xi$, where the
symplectic structure $\omega_{\xi}$ is induced from the standard
1-form of cotangent bundle $T^*N$. Then $\tilde{L}$ is a compact
Lagrangian submanifold in $SN$. There is a natural projection $p:\
SN\rightarrow N$, and each section corresponds to a splitting
$SN=N\times \R_+=N\times (e^{-\infty}, +\infty]$. The
contactomorphism $\tilde{\psi}$ can be lifted to a
$\R_+$-equivariant Hamiltonian symplectomorphism $\Psi$ on $SN$. We
denote $\Cal{L}=p^{-1}(\Lambda)$, which is also a Lagrangian
submanifold in $SN$. Then we can see that there is a 1-1
correspondence between $\tilde{L}\cap\Psi(\Cal{L})$ and
$L\cap\pi\circ\tilde{\psi}_1(\Lambda)$. However, the symplectization
$SN$ is not compact. So the ordinary method of Floer Lagrangian
intersection will be modified.

Following the argument of Ono\cite{On}, we can replace the
symplectization $(SN,\ \omega_{\xi})$ manifold by another symplectic
manifold $(Q,\ \Omega)$, which may be considered as a symplectic
filling in the negative end, so $Q$ coincides with $SN$ in the part
$N\times [e^{-C}, +\infty]\supset \tilde{L}$, where $C>0$ is a
sufficiently large number. We note that $Q$ is a 2-plane bundle over
$M$ and is diffeomorphic to the associated complex line bundle
$N\times_{S^1}\C$. We define the compatible almost complex structure
by $J'$ on $Q$ in the following way. Since $Q$ is the associated
complex line bundle, the connection $\alpha$ on $N$ gives the
decomposition of $TQ={\rm Ver}(Q)\oplus{\rm Hor}(Q)$. And we have a
$\omega$-compatible almost complex structure $J$ on $M$, then we
lift $J$ to an almost complex structure on ${\rm Hor}(Q)$. Also we
define the almost complex structure on each fiber by choosing the
standard complex structure $J_0$ on complex plane $\C$. Then we let
$J'=J\oplus J_0$, so $J'$ is uniquely determined by the choice of
$J$ on $M$ and a connection on $N$. Furthermore, Ono (c.f. section 6
in \cite{On}) showed that if we choose a generic $J$ on $M$ in the
sense of the construction of Floer homology for $(M,L)$, then $J'$
is also a regular or generic almost complex structure on $Q$. If we
write $\Pi:Q\rightarrow M$ for the natural projection, then it is a
$(J',J)$-holomorphic map. Therefore, a map
$u=\Pi\circ\tilde{u}:\Sigma\rightarrow M$ is $J$-holomorphic if and
only if $\tilde{u}:\Sigma\rightarrow M$ is $J'$-holomorphic. And we
can see that, for $r>1$, the image of the positive end $N\times
\{r\}\subset SN$ in $Q$ is $J'$-{\it convex}. So we can choose the
$\Omega$-compatible almost complex structure so that it coincides
with $J'$ outside of a compact set. For simplicity, we still denote
by $J$ this almost complex structure on $Q$ if without the danger of
confusion.

Moreover, Ono also proved that there is an a priori $C^0$-bound for
connecting orbits in $Q$ (Especially, all $J$-holomorphic curves
which we concern are contained in a compact subset $K\subset Q$,
while $K$ depends on the choice of the contact isotopy $\{\psi_t\}$
), and the bubbling off argument can go through as in the case of
compact symplectic manifold. So the minimal area of $J$-holomorphic
spheres and $J$-holomorphic discs bounding Lagrangian submanifolds
$\tilde{L}$ and $\Cal{L}$ can be achieved, we denote it by
$$\sigma(Q,\tilde{L},\Cal{L},J)=\min(\sigma_S(Q,J)|_K,
\sigma_D(Q,\tilde{L},J)|_K,\sigma_D(Q,\Cal{L},J)|_K,\sigma_D(Q,\Cal{L},\tilde{L},J)|_K)$$
and $$\sigma(Q,\tilde{L},\Cal{L})=\sup_{J_M\in {\Cal
J}}\sigma(Q,\tilde{L},\Cal{L},J=J_M\oplus J_0).$$

We will show that we can find a compactly supported Hamiltonian
diffeomorphism $\Psi'\in Ham(Q)$ such that for a compact set $K$,
the two images of $\Psi$ and $\Psi'$ coincide. For detailed
explanation, we refer to [On] or the section 2. Now we denote a
contactomorphism by $\psi$, then our main result is the following
\begin{Theorem}\label{mainthm}
Let $M$ be a compact symplectic manifold, and $N$ be the principal
$S^1$-bundle $\pi:\ N\rightarrow M$ defined above. Given a contact
isotopy $\psi_t\ |\ 0\le t\le 1\}$ on $N$, suppose $L$ is a closed
Lagrangian submanifold of $M$ which can be lifted to a Legendrian
submanifold $\Lambda$ in $N$, and
$E(\Psi')<\sigma(Q,\tilde{L},\Cal{L})$, then
$\#(L\cap\pi\circ\psi_1(\Lambda))\geq {\rm dim}H_*(L,\Z_2)$,
provided $L$ and $\pi\circ\psi_1(\Lambda)$ intersect transversally.
\end{Theorem}
\bigskip
\noindent {\bf Acknowledgement}. The author thanks Professor Yiming
Long for constant encouragement for his working on symplectic
geometry. He also thanks both referees for their careful checking
the paper and pointing out some important points and statements
which are unclearly represented in the earlier version of the paper. 
Especially, he wants to thank one of referees who pointed out one important point about possible bubbling-off of holomorphic discs from the continuation trajectories, which the author didn't take into consideration in the earlier version of the paper, and gave a suggestion of revising the paper.

\section {Preliminaries.}

We introduce some fundamental concepts and facts in symplectic and
contact geometry.

Given a $2n+1$-dimensional manifold $N$, we say $N$ is a contact
manifold if there exists a contact structure $\xi$, which is a
completely non-integrable tangent hyperplane distribution. It is
obvious that $\xi$ can locally be defined by a 1-form $\alpha$,
$i.e.\ \xi=\{\alpha=0\}$ or $\xi={\rm ker}\ \alpha$, satisfying
$\alpha\wedge(d\alpha)^n\neq 0$. If the contact structure is
co-orientable, then $\alpha$ can be global defined. We only consider
the co-oriented contact structure in this paper. The contact
manifold is denoted by $(N,\xi)$, $\alpha$ is called a contact form.
A diffeomorphism $\psi$ of $N$ is called a contactomorphism if it
preserves the co-oriented contact structure $\xi$. $\{\psi_t,\ 0\le
t\le 1\}$ is called a contact isotopy, if $\psi_0={\rm id}$ and
every $\psi_t$ is a contactomorphism. And $X_t=\frac{d\psi_t}{dt}$
is the contact vector field on $N$.

For any symplectic manifold $(M,\omega)$, there exists an almost
complex structure $J$ on $M$. We say the almost complex is {\it
compatible} with the symplectic manifold, if
$\omega(J\cdot,J\cdot)=\omega(\cdot,\cdot)$, and
$\omega(\cdot,J\cdot)>0$, which can give the Riemannian metric on
$M$.

Let $N$ be an oriented codimension 1 submanifold in an almost
complex manifold $(Q,J)$, and $\xi_x$ be the maximal $J$-invariant
subspace of the tangent space $T_xN$, then $\xi_x$ has codimension
1. And $N$ is said to be $J$-{\it convex} if for any defining 1-form
$\alpha$ for $\xi$, we have $d\alpha(v,Jv)>0$ for all non-zero
$v\in\xi_x$. This implies $\xi$ is a contact structure on $N$. It is
a fact that if $N$ is $J$-convex then no $J$-holomorphic curve in
$Q$ can touch (or tangent to) $N$ from inside (from negative side)
(c.f. \cite{G}, \cite{M}).    \\

\noindent{\bf Symplectization.}

We denote by $SN=S_{\xi}(N)$ the $\R_+$-subbundle of the cotangent
bundle $T^*N$ whose fiber at $q\in N$ are all non-zero linear forms
in $T_q^*N$ which is compatible with the contact hyperplane
$\xi_q\subset T_qN$. There is a canonical 1-form $pdq$ on $T^*N$,
and let $\alpha_{\xi}=pdq|_{SN}$, then $\omega_{\xi}=d\alpha_{\xi}$
is a symplectic structure on $SN$. Thus, we call $(SN,\
\omega_{\xi})$ the symplectization of the contact manifold $(N,\
\xi)$. We see that a contact form $\alpha:\ N\rightarrow SN$ is a
section of this $\R_+$-bundle $p:\ SN\rightarrow N$, hence we have a
splitting $SN=N\times\R_+$.

An $n$-dimensional submanifold  $\Lambda\subset(N,\xi)$ is called
Legendrian if it is tangent to the distribution $\xi$, that is to
say, $\Lambda$ is Legendrian iff $\alpha|_{\Lambda}=0$. The
preimage ${\Cal L}=p^{-1}(\Lambda)$ is an $\R_+$-invariant
Lagrangian cone in $(SN,\omega_{\xi})$. Conversely, any Lagrangian
cone in the symplectization projects onto a Legendrain submanifold
in $(N,\xi)$.

$SN$ carries a canonical conformal symplectic $\R_+$-action. Every
contactomorphism $\varphi$ uniquely lifts to a $\R_+$-equivariant
symplectomorphism $\tilde{\varphi}$ of $SN$, which is also a
Hamiltonian diffeomorphism of $SN$. Conversely, each
$\R_+$-equivariant symplectomorphism of $SN$ projects to a
contactomorphism of $(N,\xi)$. A function  $F$ on $SN$ is called a
contact Hamiltonian if it is homogeneous of degree 1, $i.e.\
F(cx)=cF(x)$ for all $c\in\R_+,\ x\in SN$.

The Hamiltonian flow generated by a contact Hamiltonian function
is $\R_+$-equivariant, it defines a contact isotopy on
$(N,\xi)$, therefore, any contact isotopy $\{\varphi_t\}$ is
generated in this sense by a uniquely defined time-dependant
contact Hamiltonian $F_t:\ SN\rightarrow \R$. There is a 1-1
correspondence between a contact vector field $X_t$ and a function
on $N$: $f_t=\alpha(X_t)$, which is also called a contact
Hamiltonian function.\\ \\
\noindent{\bf Contactization.}

If a symplectic manifold $(M,\ \omega)$ is exact, $i.e.\
\omega=d\alpha$, then it can be contactized, The contactization
$C(M,\omega)$ is the manifold $N=M\times S^1$ (or $M\times \R$)
endowed with the contact form $dz-\alpha$. Here we denote by $z$
the projection to the second factor and still denote by $\alpha$
its pull-back under the projection $N\rightarrow M$.

However, the contactization can be defined sometimes even when
$\omega$ is not exact. Suppose that the form $\omega$ represents
an integral cohomology class $[\omega]\in H^2(M)$. The
contactization $C(M,\omega)$ of $(M,\omega)$ can be constructed as
follows. Let $\pi:\ N\rightarrow M$ be a principal $S^1$-bundle
with the Euler class equal to $[\omega]$. This bundle admits a
connection whose curvature form just is $\omega$. This connection
can be viewed as a $S^1$-invariant 1-form $\alpha$ on $N$. The
non-degeneracy of $\omega$ implies that $\alpha$ is a contact form
and, therefore $\xi=\{\alpha=0\}$ is a contact structure on $N$.
The contact manifold $(N,\xi)$ is, by the definition, the
contactization $C(M,\omega)$ of the symplectic manifold
$(M,\omega)$. A change of the connection $\alpha$ leads to a
contactomorphic manifold.

We note that a Hamiltonian vector field on $(M,\omega)$ can be
lifted to a contact vector field on $N$. In fact, a Hamiltonian
function $H$ on $M$ and its Hamiltonian vector field $X_H$ satisfy
$dH=\iota(X_H)\omega$. And we know there exists a 1-1
correspondence between contact vector fields and functions on $N$,
so we obtain a contact vector field $\tilde{X}_H$ on $N$ by
$\alpha(\tilde{X}_H)=\pi^*H$. Also we have $\pi_*\tilde{X}_H=X_H$.
Thus, any Hamiltonian isotopy on $M$ is lifted to a contact
isotopy on $N$.

If $L\subset M$ is a Lagrangian submanifold, then the connection
$\alpha$ over it is flat. The pull-back $\pi^{-1}(L)\subset N$
under the projection, which is also the image of the $S^1$-action
of a Legendrian lift $\Lambda$, denoted by $\tilde{L}$, is a
Lagrangian submanifold in $SN$ and is foliated by Legendrian
leaves obtained by integrating the flat connection over $L$. If
the holonomy defined by the connection $\alpha$ is integrable over
$L$ then the Lagrangian submanifold $\tilde{L}$ is foliated by
closed Legendrian submanifolds in $N$. In particularly, this is
the case when the connection over $L$ is trivial. If this
condition is satisfied then $L$ is called exact (Bohr-Sommerfeld
condition). In this case the Lagrangian submanifold $\tilde{L}$ is
foliated by closed Legendrian lifts of $L$.

A Legendrian submanifold $\Lambda\subset (N,\xi)$ has a
neighborhood $U$ contactomorphic to the 1-jet space
$J^1(\Lambda)$. Then $\tilde{L}\cap U$ can be identified under the
contactomorphism with the so-called ``0-wall":
$W=\Lambda\times\R\subset\ J^1(\Lambda)$, which is just the set of
1-jets of function with 0 differential. \\ \\
\noindent{\bf Modify $(SN, \omega_{\xi})$.}

Now, given a contact isotopy $\{\psi_t|\ 0\le t\le 1\}$ of
$(N,\xi)$. It can be lifted to a Hamiltonian isotopy $\{\Psi_t|0\le
t\le 1\}$ of $SN$. Then, from the definition and properties listed
above, we have a 1-1 correspondence between $L\cap
\pi\circ\psi_1(\Lambda)$ and $\tilde{L}\cap\Psi_1(p^{-1}(\Lambda))$,
also they coincide with $\tilde{L}\cap\psi_1(\Lambda)$, and all
intersections are transversal. Therefore, it is natural to define
Floer homology for such a pair of Lagrangian submanifolds
$\tilde{L}$ and ${\Cal L}=p^{-1}(\Lambda)$. However, as we all know,
symplectization $SN$ is not compact, thus the ordinary method can
not directly apply. Now, we adopt Ono's argument\cite{On} to
overcome this difficulty.

We see $N$ is compact, thus there exists large $C>0$, such that the
trace of $N$ under the isotopy $\{\Psi_t|0\le t\le 1\}$ is contained
in a compact set $N\times [e^{-C},e^{C}]$, and $N\times
[e^{-C},e^{C}]$ is disjoint from $\Psi_t(SN\setminus [e^{-D},e^D]),\
t\in[0,1]$, for some number $D>C$. So the part
$N\times[e^{-D},+\infty)$ is the domain we concern. The isotopy
$\{\Psi_t\}$ is generated by a Hamiltonian $H:[0,1]\times
SN\rightarrow\R$. We can find another function $H'$, so that $H'$
equals $H$ on $N\times [e^{-C},e^{C}]$, and equals zero outside of
$N\times [e^{-D},e^{D}]$. Then we get a new Hamiltonian isotopy
$\{\Psi'_t|0\le t\le 1\}$ with compact support.

Since the boundary of the bundle $N\times[e^{-D-\epsilon},+\infty)$
is of contact type, by symplectic filling techniques, the
symplectization $(SN=N\times\R_+,\omega_{\xi})$ can be replaced by a
new symplectic manifold $(Q,\Omega)$, which is diffeomorphic to the
associated complex line bundle $N\times_{S^1}\C\rightarrow M$. In
fact, Ono showed there exists a symplectic embedding  ${\Cal F}$
from $N\times(e^{-D-\epsilon},+\infty)$ into $(Q,\Omega)$(In fact,
${\Cal F}$ is a symplectomorphism between
$N\times(e^{-D-\epsilon},+\infty)$ and
$N\times_{S^1}\C-\{0-section\}$, we refer to the appendix in
\cite{On} for details). Therefore, we just study the Lagrangian
intersection problem for $Q,\ {\Cal F}\tilde{L},\ {\Cal F}({\Cal
L}\cap N\times(e^{-D-\epsilon},+\infty))$ under Hamiltonian isotopy
$\Phi_t$ generated by a Hamiltonian defined on $Q$, which equals $\
H'\circ{\Cal F}^{-1}$ on $N\times_{S^1}\C-\{0-section\}$, and equals
zero on the 0-section. For simplicity, we still denote them by
$\tilde{L},\ {\Cal L},\ H$.

Also we notice that the positive end of $Q$ is $J$-convex, $i.e.$
for a given $E>1$, $N\times\{E\}\subset Q$ is a  $J$-convex
codimension 1 submanifold. So there is no J-holomorphic curves can
touch it, especially, there exists a $C^0$ bound for every
J-holomorphic disc $u:D^2\rightarrow Q$ with boundary in Lagrangian
submanifolds $\tilde{L}$ and $\Phi_t(\Cal{L})$(also c.f. \cite{On}).
For general case, we consider $u:\Pi=\R\times[0,1]\rightarrow Q$
with $u(\tau,0)\subset\Cal{L}$ and $u(\tau,1)\subset\tilde{L}$,
$\tau\in\R$, which is regarded as the connecting orbit between
$x_-(t)=lim_{\tau\rightarrow -\infty}u(\tau,t)$ and
$x_+(t)=lim_{\tau\rightarrow +\infty}u(\tau,t)$, solving the
perturbed Cauchy-Riemann equation $$\frac{\partial u}{\partial
\tau}=-J\frac{\partial u}{\partial t}+\nabla H(t,u(\tau,t)).$$ In
this situation, Gromov\cite{G} showed how to define an almost
complex structure $\tilde{J}_H$ on the product $\tilde{Q}=\Pi\times
Q$, such that the $\tilde{J}_H$-holomorphic sections of $\tilde{Q}$
are precisely the graph $\tilde{u}$ of solutions of the equation
above. We can see that $\tilde{Q}$ is $\tilde{J}_H$-convex, so there
is a $C^0$-bound for $\tilde{J}_H$-holomorphic curves in
$\tilde{Q}$, then the same thing happens to the connecting orbits in
$Q$.

\section {Variation and Functional.}

From the discussion above, we know that we have got a symplectic
manifold $(Q,\Omega)$, and two Lagrangian submanifolds $\tilde{L}$
and ${\Cal L}$. Then we will establish a homology theory for the
pair $(\tilde{L},\Cal{L})$ in $Q$, and study critical points of the
symplectic action functional defined on (some covering of) the space
of paths in $Q$, starting from $\Cal{L}$ with ends on $\tilde{L}$.


Let $H\in{\Cal H}(Q)$ satisfy $\|H\|<
\sigma(Q,\tilde{L},\Cal{L},J)$, and $\Psi_{(s)}^t,\ s\in [0,1]$, be
the time-$t$ flow generated by Hamiltonian $sH$ (note that
$\Psi_{(s)}^1$ is the lift of the contactomorphism $\psi_{(s)}^1$).
And denote ${\Cal L}_s=\Psi_{(s)}^1({\Cal L})$,
$\Lambda_s=\psi_{(s)}^1(\Lambda)\subset N$. We suppose that
$\tilde{L}$ intersects ${\Cal L}_1$ transversally.

Let $\Sigma$ be the connected component of constant paths in the
path space
$$\{\gamma\in C^{\infty}([0,1],Q)|\gamma(0)\in \Cal{L},
\ \gamma(1)\in\tilde{L} \}.$$ We define the closed 1-form $\alpha$
on $\Sigma$ by $$\langle
\alpha(\gamma),v\rangle=\int_0^1\Omega(\dot{\gamma}(t),v(t))dt,\
v(t)\in TQ|_{\gamma(t)},\ \forall \ t\in[0,1].$$

We also write the function $\theta:\ \Sigma\rightarrow\R$ as
$$\theta(\gamma)=-\int_0^1H(t,\gamma(t))dt.$$ Note that the zeroes
of $\alpha_s=\alpha+sd\theta$ are just time-1 trajectories generated
by the flow $\Psi_{(s)}^t$ which start from $\Cal{L}$ and end on
$\tilde{L}$. If $\gamma$ is the zero of $\alpha_s$, then the ends of
all $\gamma(1)$ are just the intersection points of $\tilde{L}$ with
${\Cal L}_s$, which are 1-1 correspondent to the zeroes of
$\alpha_s$. The purpose of this paper is to estimate from below the
number of zeroes of $\alpha_1$ .

Since $H_t$ is compactly supported on $Q$, let
$b_+=\int_0^1\max_{x\in Q}H(t,x)dt$, and $b_-=\int_0^1\min_{x\in
Q}H(t,x)dt$. Then $\|H\|=b_+-b_-$, $-b_+\le\theta(\gamma)\le-b_-$,
for all $\gamma\in\Sigma$. We introduce the Riemannian structure on
$\Sigma$ by the metric
$$(v_1,v_2)=\int_0^1\Omega(v_1(t),Jv_2(t))dt.$$ Since $$(grad_{\alpha}(\gamma),v)=\langle
\alpha(\gamma),v\rangle= \int_0^1\Omega(\dot{\gamma}(t),v(t))dt=
\int_0^1\Omega(J\dot{\gamma}(t),Jv(t))dt=(J\dot{\gamma},v),$$ so the
gradient of the closed 1-form $\alpha$ is given by $J\dot{\gamma}$,
similarly, the gradient of the closed 1-form $\alpha_s$ is
$grad_{\alpha_s}=J\dot{\gamma}-s\nabla H$.

Now, we consider the minimal covering $\pi:\
\tilde{\Sigma}\rightarrow\Sigma$ such that the form $\pi^*\alpha$ is
exact, $i.e.$ there is a functional $F$ on $\tilde{\Sigma}$, such
that $\ \pi^*\alpha=dF$, and its structure group $\Gamma$ is free
abelian. Denote $F_s=F+s(\theta\circ\pi)$, so $dF_s=\pi^*\alpha_s$.
The gradient $\nabla F_s$ of the functional $F_s$, with respect to
the lift of the Riemannian structure on $\Sigma$, is a
$\Gamma$-invariant vector field on $\tilde{\Sigma}$, and
$\pi_*\nabla F_s=grad_{\alpha_s}$. Then we consider the moduli space
of thus gradient flows connecting a pair of critical points
$(x_-,x_+)$ of $F_s$
$$M_s(x_-,x_+)=\{ u:\ \R\rightarrow\tilde{\Sigma}|
\ \frac{du(\tau)}{d\tau}=-\nabla F_s(u(\tau)),\ u\ {\rm is\ not\
constant},\  \lim_{\tau\rightarrow\pm\infty}u(\tau)=x_{\pm} \}.$$

Denote by ${\Cal M}_s=\bigcup_{x_{\pm}}M_s(x_-,x_+)$ the collection,
and the nonparameterized space by
$\hat{M}_s(x_-,x_+)=M_s(x_-,x_+)/\R$, and the natural quotient map
$q:\ M_s\rightarrow\hat{M}_s $. Choosing a regular
$\Omega$-compatible almost complex structure $J$ on $Q$ (c.f.
\cite{On})\footnote{Recall that the $J$ used here is just the
$J'=J\oplus J_0$ given in the Introduction part, by generic
choosing $\omega$-compatible almost complex structure $J$ on $M$ we
can obtain the regular or generic $\Omega$-compatible structure $J'$
on $Q$. The arguments in \cite{On} for $J'$-holomorphic maps can
apply to our $H$-perturbed $J'$-holomorphic map by similar
statements as those in \cite{FHS}. We can overcome the similar
problem which appears in the continuation argument of Section 6.}, we may
assume that there is a dense set $T\subset[0,1]$ such that for all
$s\in T$, $M_s(x_-,x_+)$ are finite dimensional smooth manifolds,
consequently, $\tilde{L}$ intersects ${\Cal L}_s$ transversally.

We define the length of a gradient trajectory $u\in M_s(x_-,x_+)$
by $l_s(u)=F_s(x_-)-F_s(x_+)$. If $\hat{u}\in \hat{M}_s$, then we
define its length naturally by $l_s(\hat{u})=l_s(u)$, where
$\hat{u}=q\circ u$. Denote $\Pi=\R\times[0,1]$, then the map
$\bar{u}:\ \Pi\rightarrow Q$, defined by
$\bar{u}(\tau,t)=\pi(u(\tau))(t)$, satisfies the following
perturbed Cauchy-Riemann equation
$$\frac{\partial\bar{u}(\tau,t)}{\partial\tau}=
-J(\bar{u}(\tau,t))\frac{\partial\bar{u}(\tau,t)}{\partial t}+
s\nabla H(t,\bar{u}(\tau,t)),$$ with limits
$$\lim_{\tau\rightarrow\pm\infty}\bar{u}(\tau,t)=\pi(x_\pm)=\bar{x}_\pm(t).$$
It is easy to see that
$l_0(u)=\int_{-\infty}^{+\infty}u^*dF=\int_{\Pi}\bar{u}^*\Omega$.

If $u\in M_0$, then $\bar{u}$ is a J-holomorphic map from $\Pi$ to
$Q$. From Oh's removing of boundary singularities theorem\cite{Oh1},
$\bar{u}$ can be extended to a J-holomorphic curve $\bar{u}':\
(D^2,\partial^+ D^2,\partial^-
D^2)\rightarrow(Q,\tilde{L},\Cal{L})$, where $D^2=\bar{\Pi}$ is the
two-point compactification of $\Pi$. Since
$l_0(u)=\int_{\Pi}\bar{u}^*\Omega=\int_{D^2}(\bar{u}')^*\Omega$, we
know that $l_0(u)\geq \sigma_D(Q,\tilde{L},\Cal{L},J)$.

\section {Define and Compute Homology for $C_{\varepsilon}^0$}

We denote by $Y_s$ the set of critical points of $F_s$, and by
$C_s$ the vector space spanned by $Y_s$ over $\Z_2$.

Since $Y_s$ is $\Gamma$-invariant, $C_s$ has a structure of free
$K$-module with rank$=\#(\tilde{L}\cap{\Cal L}_s)$, $s\in T$, where
$K=\Z_2[\Gamma]$. Our aim in this section is to establish some
homology for the complex $C_{\varepsilon}$, where $\varepsilon$ is
small enough. We write the following definition similar as the one
given by Chekanov\cite{C2}.
\begin{Definition}\label{short}
Fix $\delta>0$, satisfying
$\Delta:=\|H\|+\delta<\sigma(Q,\tilde{L},\Cal{L},J)$. A gradient
trajectory $u\in M_s$ is said to be short if $l_s(u)\le\Delta$, and
be very short if $l_s(u)\le\delta$.
\end{Definition}

Now we denote the area by $A(u)=\int_{\Pi}\bar{u}^*\Omega$, and
$h(u)=s\int_{-\infty}^{+\infty}u^*d(\theta\circ\pi)$, then still
write $l(u)=l_s(u)=A(u)+h(u)$, we have
\begin{Lemma}\label{4.1}
If $u$ is very short $i.e.$ $l(u)\le \delta$, then the area
$A(u)\le \Delta$.
\end{Lemma}
Proof. since $\theta=-\int_0^1H(t,\gamma(t))dt\in [-b_+,-b_-]$,
then
$$h(u)=s\int_{-\infty}^{+\infty}u^*d(\theta\circ\pi)
=s\theta(\pi(u(\tau)))|_{-\infty}^{+\infty}\geq s(b_--b_+),$$ so
$$A(u)=l(u)-h(u)\le \delta-(b_--b_+)=\|H\|+\delta=\Delta. $$

\noindent Then we can prove the following key
lemma\footnote{Actually, the lemma is essentially proved by Chekanov
(c.f. Lemma 6 in \cite{C2}), here we rewrite it in our settings with
some modifications.}
\begin{Lemma}\label{U}.
For a small neighborhood $U$ of $\tilde{L}$ in $Q$, there exists a
$\varepsilon_0>0$, such that for any positive
$\varepsilon<\varepsilon_0$, every short gradient trajectory $u\in
M_{\varepsilon}$ is very short, and for every short $u$ we have
$\bar{u}(\Pi)\subset U$.
\end{Lemma}
Proof. We prove it by contradiction. For the first claim, we suppose
there is a sequence $u_n\in M_{s_n}$ and a positive number $c$ with
$\delta\le c\le \Delta$ so that when $s_n\rightarrow 0$ then
$l_{s_n}(u_n)\rightarrow c$. By Gromov's compactness theorem, there
are some subsequence of $\bar{u}_n=\pi(u_n)$ convergent to
$\bar{u}_\infty$ which is a collection of $J$-holomorphic spheres
and $J$-holomorphic discs bounding $\tilde{L}$ and/or $\Cal{L}$.
Then the total symplectic area of this limit collection is just
$l_0(u_\infty)=c$ which by the assumption of Theorem 1 is lager than
$\sigma(Q,\tilde{L},\Cal{L})$, but $c\le
\Delta<\sigma(Q,\tilde{L},\Cal{L})$, so the claim holds. For the
second claim, the argument is similar. Note that if the image
$\bar{u}_\infty(\Pi)$ of the limit collection is not contained in
$U$, then at least one of the $J$-curve is not contained in $U$
which is nonconstant and its area will be larger than
$\sigma(Q,\tilde{L},\Cal{L},J)>\Delta$, this contradicts the (very)
shortness condition.\ \ \ \ QED.

\bigskip
Then, we denote by $M'_{\varepsilon}\
(\hat{M}'_{\varepsilon})\subset M_{\varepsilon}\
(\hat{M}_{\varepsilon})$ the set of all short gradient trajectories
(nonparameterized short gradient trajectories). And we can define
the $\Z_2$-linear map $\partial:\ C_{\varepsilon}\rightarrow
C_{\varepsilon}$ by
$$\partial(x)=\sum_{y\in Y_{\varepsilon}}\#\{{\rm isolated\ points
\ of}\  \hat{M}'_{\varepsilon}(x,y) \}y,$$ for $\forall\ x\in
Y_{\varepsilon}$.

Let $\varepsilon\in T$ be sufficiently small and satisfy the
conditions of lemma \ref{U}. Choose an element $x_0\in
Y_{\varepsilon}$, then we can define a subclass
$Y_{\varepsilon}^0\subset Y_{\varepsilon}$ by
$$Y_{\varepsilon}^0=\{x\in Y_{\varepsilon}\ |\ |F_{\varepsilon}
(x)-F_{\varepsilon}(x_0)|\le \delta \}.$$ Then we see that the
projection $\pi$ bijectively maps $Y_{\varepsilon}^0$ onto the set
of zeroes of the form $\alpha_{\varepsilon}$. And we get the
bijection
$$Y_{\varepsilon}^0\times\Gamma\rightarrow Y_{\varepsilon}:
\ \ (y,a)\mapsto a(y),$$ which induces the isomorphism
$C_{\varepsilon}^0\otimes K\rightarrow C_{\varepsilon}$, where
$C_{\varepsilon}^0\subset C_{\varepsilon}$ is spanned over $\Z_2$
by $Y_{\varepsilon}^0$. Now, for sufficiently small
$\varepsilon\in T$, we can establish the homology for
$(C_{\varepsilon},\partial)$

\begin{Lemma}\label{homology}
$1^{\circ}$ The map $\partial$ is $K$-linear, well defined, and
$\partial(C_{\varepsilon}^0)\subset C_{\varepsilon}^0$;\\
$2^{\circ}$ If $\varepsilon\in T$ is sufficiently small, then
$\partial^2=0$;\\
$3^{\circ}$ The homology $H(C_{\varepsilon}^0,\partial)\cong
H_*(\Lambda,\Z_2)$.
\end{Lemma}
Proof. $1^{\circ}$ Since the gradient flow is $\Gamma$-invariant,
$\partial$ is naturally $K$-linear. We know that the bubbling off
can not occur. Indeed, since $\varepsilon$ is sufficiently small,
then $u\in M_{\varepsilon}$ is very short, $l_{\varepsilon}(u)\le
\delta$, by the lemma \ref{4.1}, the area
$A(u)\le\Delta<\sigma(Q,\tilde{L},{\Cal L},J)$, and from the
assumption in our theorem, the area of any J-holomorphic sphere or
J-holomorphic disc bounding $\tilde{L}$ and ${\Cal L}$ is larger
than $\sigma(Q,\tilde{L},{\Cal L},J)$. Thus,
$\hat{M}'_{\varepsilon}(x,y)$ is compact
and the number of its isolated points is finite.\smallskip\smallskip\\
$2^{\circ}$ Suppose $\varepsilon\in T$ satisfy the conditions in
lemma \ref{U}. If $\|H\|=0$, $\Delta=\delta$, then $H\equiv
const.$ and $\psi_H\equiv id$, it is a trivial case. If $\|H\|>0$,
we can always choose a fixed
$\delta<\frac{1}{2}\|H\|<\frac{1}{2}\Delta$. Consider a pair of
isolated trajectories $u_1\in \hat{M}'_{\varepsilon}(x,y)$,
$u_2\in \hat{M}'_{\varepsilon}(y,z)$. Then there exists a unique
1-dimensional connected component ${\Cal C}\subset
\hat{M}_{\varepsilon}(x,z)$ such that $(u_1,u_2)$ is one of the
two ends of compactification of ${\Cal C}$(c.f. \cite{F1}). Since
the length is additive under gluing, we have, for $\forall
u\in{\Cal C}$, $l_{\varepsilon}(u)=l_{\varepsilon}(u_1)+
l_{\varepsilon}(u_2)<2\delta<\Delta$. By lemma \ref{U}, $u$ is
also very short. From the lemma \ref{4.1}, we know the bubbling
off doesn't occur, too. Then the other end of ${\Cal C}$ can be
compactified by a pair of isolated trajectories
$u'_1\in\hat{M}_{\varepsilon}(x,y)$,
$u'_2\in\hat{M}_{\varepsilon}(y,z)$. Also $l_{\varepsilon}(u'_1)+
l_{\varepsilon}(u'_2)=l_{\varepsilon}(u_1)+
l_{\varepsilon}(u_2)<\Delta$, thus $u'_1,\ u'_2$ are short
trajectories. So we know that, for $x,z\in Y_{\varepsilon}$, the
number of isolated trajectories $(u_1,u_2)\in
\hat{M}'_{\varepsilon}(x,y)\times\hat{M}'_{\varepsilon}(y,z)$ is
even, $i.e.\ \partial^2(x)=0$.\smallskip\smallskip\\
$3^{\circ}$ From the lemma \ref{U}, we know that for any
$u\in\hat{M}'_{\varepsilon}(x,y)$, $\bar{u}(\bar{\Pi})\subset U$.
That is to say, if $\varepsilon$ is small enough, then
$\Lambda_{\varepsilon}=\psi_{(\varepsilon)}^1(\Lambda)$ is always
contained in a small neighborhood $U'=U\cap N$ of the Legendrian
submanifold $\Lambda$ in $N$. By Darboux's theorem, $U'$ is
contactomorphic to a neighborhood of the 0-section in the 1-jet
space $J^1(\Lambda)$. This contactomorphism moves $\tilde{L}\cap U'$
onto the``0-wall" $W$, $i.e.$ the space of 1-jets of functions with
0 differential. Thus a Legendrian submanifold $\Lambda'$, which is
$C^1$-close to $\Lambda$ and transverse to $W$, corresponds to a
Morse function $\beta:\ \Lambda\rightarrow\R$ so that the
intersection points of $\tilde{L}$ and $\Lambda'$ are in 1-1
correspondence with the critical points of the function $\beta$. We
can explicitly choose a metric on $\Lambda$ and a generic almost
complex structure $J$ (recall the footnote in section 3) on the
symplectization $SN$ in such a way that the gradient trajectories in
${\Cal M}_{\varepsilon}$ would be in 1-1 correspondence with the
gradient trajectories of the function $\beta$ connecting the
corresponding critical points of this function. Thus we can identify
our complex $C^0_{\varepsilon}$ with the Morse chain complex for the
function $\beta$ (here  we may also reduce the problem to Lagrangian
intersections in $M$ by applying continuation argument and
projecting the manifold to $M$, the method of equating Floer and
Morse complex is standard, we refer the reader to \cite{S}), so we
have an isomorphism $H(C^0_{\varepsilon},\partial)\cong
H_*(\Lambda,\Z_2)$.

\section {Homology Algebra.}

Under the condition $\pi_2(Q,\cdot)=0$ or the monotonicity
assumption\cite{F1}\cite{Oh2}, the Floer homology of the complex
$C_s,\ s\in T\subset[0,1]$, can be defined, $i.e.$
$HF_*(C_s,\partial)$. Then we can use the classical continuation
method (c.f. \cite{F4}\cite{M90}\cite{Oh2}) to prove the isomorphism
between $HF_*(C_{\varepsilon},\partial)$ and $HF_*(C_1,\partial)$,
that means to construct chain homotopy $\Phi'\Phi\sim
id_{\varepsilon}$ (and $\Phi\Phi'\sim id_1$), where $\Phi:\
(C_{\varepsilon},\partial)\rightarrow(C_{1},\partial)$, $\Phi':
(C_{1},\partial)\rightarrow(C_{\varepsilon},\partial)$ are chain
homomorphisms defined similarly as the definition of $\partial$,
except for considering the moduli space of continuation
trajectories. That is to say, in order to prove $\Phi'\Phi\sim
id_{\varepsilon}$, we should show there exists a chain homomorphism
${\bf h}:\
(C_{\varepsilon},\partial)\rightarrow(C_{\varepsilon},\partial)$, so
that $$\Phi'\Phi-id={\bf h}\partial-\partial {\bf h}.$$

However, in general case, we can not define appropriately any
homology for $C_s$ unless $s$ is small enough. Then we may only
prove a weaker ``homotopy", which is called $\lambda$-{\it
homotopy} by Chekanov. In fact, for the aim of estimating from
below the number of critical points of the functional $F_s$, this
$\lambda$-homotopy is enough.

We shall use the following homology algebraic result introduced
and proved by Chekanov \cite{C2}.

Let $\Gamma$ be a free abelian  group equipped with a monomorphism
$\lambda:\ \Gamma\rightarrow \R$, which we call a {\it weight
function}. Denote
$$\Gamma^+=\{a\in\Gamma|\lambda(a)>0\}\ \ \ \
\Gamma^-=\{a\in\Gamma|\lambda(a)<0\}$$ Let $k$ be a communicative
ring. Consider the group ring $K=k[\Gamma]$. For a $k$-module $M$,
we have the natural decomposition $M\otimes K=M^+\oplus M^0\oplus
M^-$. where $M^+=\Gamma^+(M)$, $M^0=M$, $M^-=\Gamma^-(M)$.
Consider the projections $$p^+:\ M\otimes K\rightarrow M^+\oplus
M^0,\ \ \ \ p^-:\ M\otimes K\rightarrow M^0\oplus M^-.$$ Assume
that $(M,\partial)$ is a differential $k$-module, then $\partial$
naturally extends to a $K$-linear differential on $M\otimes K$.
\begin{Definition}\label{homotopy}
We say two linear maps $\phi_0,\ \phi_1:\ M\otimes K\rightarrow
M\otimes K$ are $\lambda$-homotopic if there exists a $K$-linear
map ${\bf h}:\ M\otimes K\rightarrow M\otimes K$ such that
$$p^+(\phi_0-\phi_1+{\bf h}\partial+\partial {\bf h})p^-=0.$$
\end{Definition}

\begin{Lemma}\label{Ch}{\rm \cite{C2}}
Let $\lambda$ be a weight function on a free abelian group
$\Gamma$. Assume $(M,\partial)$ to be a differential $k$-module
and $N$ to be a $K$-module, where $K=k[\Gamma]$. if the maps
$\Phi^+:\ M\otimes K\rightarrow N,\ \Phi^-:\ N\rightarrow M\otimes
K$ are $K$-linear and $\Phi^-\Phi^+$ is $\lambda$-homotopic to the
identity, then ${\rm rank}_KN\geq {\rm rank}_kH(M,\partial)$.
\end{Lemma}

\section {Proof of Theorem 1.}

Given a $(s_-,s_+)$ continuation function $\rho:\R\rightarrow[0,1]$
satisfying
$$\rho(\tau)=\left\{
\begin{array}{cc}
s_-,& \ {\rm if}\ \tau<-r,\\
s_+,&  {\rm if}\ \tau>r,
\end{array}
\right.
$$
where $r\in\R$, we can define the moduli space of
continuation trajectories
$$M_{\rho}(x_-,x_+)=\{u:\R\rightarrow\tilde{\Sigma}\ |
\ \frac{du(\tau)}{d\tau}=-\nabla F_{\rho(\tau)}(u(\tau)), \
\lim_{\tau\rightarrow\pm\infty}u(\tau)=x_{\pm}\},$$ where
$x_{\pm}\in Y_{s_{\pm}}$. And we denote the collection by ${\Cal
M}_{\rho}=\bigcup_{x_-,x_+} M_{\rho}(x_-,x_+)$. The length of
a continuation trajectory is defined by
$l_{\rho}(u)=F_{s_-}(x_-)-F_{s_+}(x_+)$.

Choose a monotone $(\varepsilon,1)$ continuation function $\rho_+$
and a monotone $(1,\varepsilon)$ continuation function $\rho_-$.
For generic $H$, $M_{\rho_{\pm}}(x_-,x_+)$ are smooth manifolds.
We will say a continuation trajectory $u\in
M_{\rho_{+}}(x_-,x_+)$ (or $u\in
M_{\rho_{-}}(x_-,x_+)$) is short if $l_{\rho_{+}}(u)\le\delta+(1-\varepsilon)b_+$ (resp. $l_{\rho_{-}}(u)\le\delta+(\varepsilon-1)b_-$).
The subspace of all short trajectories is denoted by
$M'_{\rho_{\pm}}(x_-,x_+)\subset {\Cal M}_{\rho_{\pm}}$.

Then, under an ideal assumption\footnote{For general case, we have to modify the continuation map. The possibility that there exists a sequence of continuation trajectories reaching the negative end was pointed out to the author by one of referees.}, $i.e.$ if there is no a sequence of continuation trajectories reaches the negative
end ($i.e$ the zero section of $Q\rightarrow M$), we can simply construct the $\Z_2$-linear continuation map
$\Phi^+:C_{\varepsilon}\rightarrow C_1,\ \Phi^-:C_1\rightarrow
C_{\varepsilon}$ as
$$\Phi^+(x)=\sum_{y\in Y_1}\#\{{\rm isolated\
points\ of\ } M'_{\rho_{+}}(x,y)\}y,$$
$$\Phi^-(y)=\sum_{z\in Y_{\varepsilon}}\#\{{\rm isolated\
points\ of\ } M'_{\rho_{-}}(x,z)\}z,$$ where $x\in
Y_{\varepsilon}$. The following lemma implies that the definition
of $\Phi^{\pm}$ above is sound.

\begin{Lemma}\label{finite}
If $u\in M_{\rho_{+}}(x_-,x_+)$, then $l_{\rho_{+}}(u)\geq
(1-\varepsilon)b_-$. If $u\in M_{\rho_{-}}(x_-,x_+)$, then $l_{\rho_{-}}(u)\geq
(\varepsilon-1)b_+$. And the sum in the definition of $\Phi^{\pm}$ is
finite.
\end{Lemma}
Proof. \ Recall $F_s=F+s\theta\circ\pi$, $s\in[0,1]$. For a
$(s_-,s_+)$ continuation function $\rho:\R\rightarrow[0,1]$, we
have $F_{\rho(\tau)}=F+\rho(\tau)\theta\circ\pi$. So the length of
a continuation trajectory is
$$l_{\rho}(u)=F_{s_-}(x_-)-F_{s_+}(x_+)=
-\int_{-\infty}^{+\infty}u^*dF_{\rho(\tau)}=
-\int_{-\infty}^{+\infty}u^*dF-\int_{-\infty}^{+\infty}
u^*d(\rho(\tau)\theta\circ\pi)$$
$$=A(u)+h(u),$$
 where we denote $$A(u)=-\int_{-\infty}^{+\infty}u^*dF=
 \int_{\Pi}\bar{u}^*\Omega=\int_{-\infty}^{+\infty}
 (\frac{du(\tau)}{d\tau},\frac{du(\tau)}{d\tau})d\tau=
 \int_{-\infty}^{+\infty}\|\nabla F_{\rho(\tau)}\|^2d\tau\geq 0,$$
 and $$h(u)=-\int_{-\infty}^{+\infty}u^*d(\rho(\tau)\theta\circ\pi)=
-\int_{-\infty}^{+\infty}
\frac{d\varrho(\tau)}{d\tau}\theta(\pi(u(\tau)))d\tau.$$ Recall
that $$\theta=-\int_{0}^{1}Hdt\in[-b_+,-b_-],$$ if $u\in
M_{\rho_+}(x_-,x_+)$, $$l(u)\geq h(u)\geq (1-\varepsilon)b_-;$$ if $u\in
M_{\rho_-}(x_-,x_+)$, $$l(u)\geq h(u)\geq (\varepsilon-1)b_+.$$
 Thus, For a short trajectory $u\in
 M'_{\rho_{\pm}}(x_-,x_+)$, $$A(u)=l(u)-h(u)\le \delta+(1-\varepsilon)(b_+-b_-)
 =\|H\|+\delta=\Delta.$$ Since $A(u)=
 \int_{\Pi}\bar{u}^*\Omega\le\Delta<\sigma(Q,\tilde{L},J)$,
by Gromov's arguments, no bubbling can occur, then spaces
$M'_{\rho_{\pm}}(x_-,x_+)$ are compact, and the finiteness of
$\#\{{\rm isolated\ points\ of\ }
M'_{\rho_{\pm}}(x,y)\}$ is verified.

However, in general case, the ideal assumption is not always satisfied. If a sequence of continuation trajectories converges to a curve reaching the negative end of $\Cal L$, the boundary of moduli space of continuation trajectories will contain such curve. So we need modify the definition of $\Phi^{\pm}$ to get a well-defined continuation map. In fact, Ono (c.f. \cite{On} P.218) had dealt with a similar problem by considering the algebraic intersection number of the continuation trajectories with the zero section $O_M$ of $Q\rightarrow M$. Under Ono's assumption $\pi_2(M,L)=0$, bubbling off of holomorphic discs contained in the zero section of $Q$ with boundary on $L$ never occurs. In our case, since we have the bound for energy, such kind of bubbling off of discs is also avoided. Thus, we can define homomorphisms $\Phi_k^+:C_{\varepsilon}\rightarrow C_1,\ \Phi_k^-:C_1\rightarrow C_{\varepsilon}$ as 
$$\Phi_k^+(x)=\sum_{y\in Y_1}Int^+_{2,k}(x,y)y,$$
$$\Phi_k^-(y)=\sum_{z\in Y_{\varepsilon}}Int^-_{2,k}(y,z)z,$$
where $Int^{+}_{2,k}(x,y)$ ($Int^-_{2,k}(y,z)$) is the mod-2 number of the isolated points $u$ in the continuational moduli space $M'_{\rho_{+}}(x,y)$ (resp. $M'_{\rho_{-}}(y,z)$) which have the algebraic intersection number $u\cdot O_M=\frac{k}{2}$.

With the restriction of bound of energy, we can make similar discussions as in \cite{On} to get the finiteness of $Int^{\pm}_{2,0}$. So $\Phi^{\pm}_0$ are just our favorite continuation maps. We will not list the detailed arguments here and refer the reader to \cite{On} for original discussion. In the following, we will still denote the continuation maps by $\Phi^{\pm}$ for simplicity.\\

Then we can use the homology algebraic result listed in the
section 5 to prove the theorem 1, provided there exists a
$\lambda$-homotopy between $\Phi^-\Phi^+$ and the identity.\\

\noindent$\bullet$ {\it Prove the theorem 1}.

Now in our case, let $k=\Z_2,\ K=\Z_2[\Gamma],\
M=C_{\varepsilon}^{0},\ M\otimes K=C_{\varepsilon},\ N=C_1$, and
$\Gamma$ be the structure group of the covering. The weight
function $\lambda:\Gamma\rightarrow\R$ can be defined as
$\lambda(a)=F(a(x))-F(x)$. We also have decompositions
$Y_{\varepsilon}=Y_{\varepsilon}^+\cup Y_{\varepsilon}^0\cup
Y_{\varepsilon}^-$, $C_{\varepsilon}=C_{\varepsilon}^+\oplus
C_{\varepsilon}^0\oplus C_{\varepsilon}^-$, where
$Y_{\varepsilon}^{\pm}=\Gamma^{\pm}(Y_{\varepsilon}^0)$.

Assume that we have got a $\lambda$-homotopy
${\bf h}:C_{\varepsilon}\rightarrow C_{\varepsilon}$, then by lemma
\ref{Ch} and lemma \ref{finite} we have
$$\#(L\cap\pi\circ\psi_1(\Lambda))=\#(\tilde{L}\cap\psi_1(\Lambda))=\#(\tilde{L}\cap\Psi_1({\Cal L}))=
{\rm rank}_KC_1$$
$$\geq {\rm
rank}_kH(C_{\varepsilon}^0,\partial)={\rm dim}H_*(\Lambda,\Z_2)={\rm
dim}H_*(L,\Z_2).\ \ \ \ \ \ \ \ \ \ \ \ \ \ \ $$
This finishes the proof of the Theorem 1.\\

In the rest of this section, we show a sketchy proof of the
existence of $\lambda$-homotopy.
\begin{Lemma}\label{last}
There exists a $\lambda$-homotopy ${\bf h}:C_{\varepsilon}\rightarrow
C_{\varepsilon}$ between $\Phi^-\Phi^+$ and the identity.
\end{Lemma}
Proof. We follow the arguments of Chekanov and state his main
thought. For constructing the homomorphism ${\bf h}$, we use a family of
$(\varepsilon,\varepsilon)$ continuation functions $\mu_c,\
c\in[0,+\infty)$ satisfying

\noindent 1) $\mu_0(\tau)\equiv\varepsilon$,

\noindent 2) $\frac{du_c(\tau)}{d\tau}\geq 0$, if $\tau<0$;
$\frac{du_c(\tau)}{d\tau}\le 0$, if $\tau>0$,

\noindent 3) $c\mapsto \mu_c(0)$ is a monotone map from
$[0,+\infty)$ onto $[\varepsilon,1]$,

\noindent 4) when $c$ is large enough $\mu_c(\tau)=\left\{
\begin{array}{cc}
\rho_+(\tau+c),& {\rm if}\ \tau\le 0;\\
\rho_-(\tau-c),& {\rm if}\ \tau\geq 0.
\end{array}
\right. $


Then we denote the moduli space
$$M_{\mu}(x_-,x_+)=\{(c,u)|u\in M_{\mu_c}(x_-,x_+)\},
\ x_{\pm}\in Y_{\varepsilon}.$$ For generic $H$,
$M_{\mu}(x_-,x_+)$ are smooth manifolds. 

And like the arguments for $(\varepsilon,1)$ or $(1,\varepsilon)$-continuation trajectories shown before in this section, under a similar ideal assumption, $i.e.$ no any sequence of $\mu_c$-continuation trajectories reaches the negative end of $\Cal L$, we can define the $\Z_2$-linear map for $C_{\varepsilon}^0$ as\footnote{Otherwise, we will again adopt the Ono's argument to take into consideration of the algebraic intersection number, the way of modifying the definition of the map ${\bf h}$ is similar as the way of modifying $\Phi^{\pm}$ we have stated before. For simplicity, we just show the argument under this ideal assumption.} 
$${\bf h}(x)=\sum_{y\in Y_{\varepsilon}^0}\#
\{{\rm isolated\ points\ of\ }M'_{\mu}(x,y)\}y,\ \ x\in
Y_{\varepsilon}^0$$ where $M'_{\mu}(x,y)$ is the subset of the moduli space $M_{\mu}(x,y)$ which contains only short 
$\mu_c$-continuation trajectories, a $\mu_c$-continuation trajectory $u\in
M_{\mu_c}(x_-,x_+)$ is called short if its length $l(u)=l_{\mu_c}(u)\le \delta$.
Moreover, For any $u\in M_{\mu_c}(x_-,x_+)$, we have $l_{\mu_c}(u)=A(u)+h(u)\geq h(u)\geq (\mu_c(0)-\varepsilon)b_-+(\varepsilon-\mu_c(0))b_+\geq b_--b_+= -\|H\|$, and if
$l(u)\le\delta$, then $A(u)=l(u)-h(u)\le\delta+\|H\|\le \Delta$.
  
The map ${\bf h}$ can be extended naturally to a $K$-linear
map on $C_{\varepsilon}$. Since for $u\in M'_{\mu}(x,y)$, $l_{\mu_c}(u)\le\delta,\
A(u)\le\Delta$, the bubbling off does not occur, $M'_{\mu}(x,y)$
is compact and the sum is finite, thus the map ${\bf h}$ is well defined. 

To prove ${\bf h}$ is a $\lambda$-homotopy, we have to verify
$$p^+(x+\Phi^-\Phi^+x+{\bf h}\partial x+\partial {\bf h}x)=0,$$ for
$\forall\ x\in Y_{\varepsilon}^0\cup Y_{\varepsilon}^-$. This will
follow from the standard gluing argument involving the ends of the
1-dimensional part $\aleph$ of $M_{\mu}(x,z),\ z\in
Y_{\varepsilon}^+\cup Y_{\varepsilon}^0$. Since $x\in
Y_{\varepsilon}^0\cup Y_{\varepsilon}^-,\ z\in
Y_{\varepsilon}^+\cup Y_{\varepsilon}^0$, and $\varepsilon$ is
small enough, we know $l(u)\le \delta$ for $u\in\aleph$. Indeed,
$l(u)=F_{\varepsilon}(x)-F_{\varepsilon}(z)$, and there exist $x'$
and $z'$ in $Y_{\varepsilon}^0$ and $a\in \Gamma^+\cup\Gamma^0,\
b\in\Gamma^0\cup\Gamma^-$ such that $z=a(z'),\ x=b(x')$, and
$F_{\varepsilon}(x)-F_{\varepsilon}(x')=\lambda(b)\le 0$,
$F_{\varepsilon}(z')-F_{\varepsilon}(z)=-\lambda(a)\le 0$, also we
know that $F_{\varepsilon}(x')-F_{\varepsilon}(z')\le \delta$
since $x',z'\in Y_{\varepsilon}^0$, this implies $l(u)\le \delta$,
so $A(u)\le\Delta$\footnote{From the proof of $l(u)\le \delta$
here, the reader can see why we would only verify the so called
$\lambda$-homotopy. }. This disappears the bubbling off. Then the
compactification of $\aleph$ shows that the left hand side of the
formula above has the expression
$$\sum_{z\in Y_{\varepsilon}^+\cup Y_{\varepsilon}^0}\# \{S(x,z)\}z,$$
and the number $\# \{S(x,z)\}$ is even (the one more thing we
should verify than the standard gluing argument is to prove the
other ends of the compactification which are pairs of continuation
trajectories are still all short, this is not difficult to do\footnote{For the modified continuation maps $\Phi^{\pm}_0$, we can also verify the continuation trajectories in the other ends have the same 0-algebraic intersection number with the zero section $O_M$, see \cite{On}.}).
This ends the proof of the lemma and the theorem 1. For more
details, the reader may refer  to \cite{C2}\cite{F4}\cite{M90}.\\

\bigskip
\noindent\small Institute of Mathematical Science\\
\noindent\small Nanjing University\\
\noindent\small Nanjing 210093, P.R.China\footnote{This is the
present postal address of the author.}\\
\noindent\small And\\
\noindent\small Chern Institute of Mathematics\\
\noindent\small Nankai University\\
\noindent\small Tianjin 300071, P.R.China\\
\noindent\small E-mail:  hailongher@ims.nju.edu.cn

\end{document}